\documentclass[11pt]{article}
\usepackage{amsmath}

\makeatletter
\renewcommand\footnoterule{%
  \kern-3\p@
  \hrule\@width \textwidth
  \kern2.6\p@}

          \usepackage{amsmath}
          \usepackage{amsthm}
\date{}          

\begin{document}

\title{\bf{On the Logarithmic Coefficients  of Close to Convex Functions}}

\author{\bf{D K Thomas$^\ast$}}
 
\maketitle

\footnote{$^\ast$Department of Mathematics, Swansea University, Singleton Park, Swansea, SA2 8PP, UK. d.k.thomas@swansea.ac.uk}

\begin{abstract}
For $f$ analytic  and close to convex in $D=\{z: |z|< 1\}$, we give sharp estimates for the logarithmic coefficients $\gamma_{n}$ of $f$ defined by  $\log \dfrac{f(z)}{z}=2\sum_{n=1}^{\infty} \gamma_{n}z^{n}$ when $n=1, 2,3$.
\end{abstract}
\bigskip


\bigskip

\noindent {\bf Introduction}

\bigskip

Let S be the class of normalised analytic univalent functions $f$ for $z\in D=\{z:|z|<1\}$ and given by

\begin{equation*}
f(z)=z+\sum_{n=2}^{\infty}a_{n}z^{n}.
\end{equation*}

The logarithmic coefficients of $f$ are defined in $D$ by

\begin{equation}
\log \frac{f(z)}{z}=2\sum_{n=1}^{\infty} \gamma_{n}z^{n}.
\end{equation}
The logarithmic coefficients $\gamma_{n}$ play a central role in the theory of univalent functions. Milin conjectured that for $f\in S$ and $n\ge2$,

\begin{equation}
\sum_{m=1}^{n}\sum_{k=1}^{m}(k|\gamma_{k}|^2-\frac{1}{k})\le 0
\end{equation}

\noindent and it is not difficult to see that $(2)$ implies the Bieberbach conjecture.  It was a proof of $(2)$ that De Branges established in order to prove the conjecture. 
\bigskip

Very few exact upper bounds for $\gamma_{n}$ seem have been established, with more attention being given to results of an average sense (see e.g. [1, 2]). Moreover it is known that for $f\in S$, the expected inequality $|\gamma_{n}|\le \dfrac{1}{n}$ is false even in order of magnitude [1, Theorem 8.4]. 

\bigskip

\noindent Differentiating $(1)$ and equating coefficients gives

\begin{gather}
\gamma _{1}=\frac{1}{2}a_{2}\\
\gamma _{2}=\frac{1}{2}(a_{3}-\frac{1}{2}a_{2}^{2})\\
\gamma _{3}=\frac{1}{2}(a_{4}-a_{2}a_{3}+\frac{1}{3}a_{2}^{3})
\end{gather}

\noindent Hence $|\gamma_{1}|\le 1$  follows at once from $(3)$, and use of the Fekete-Szeg\"o inequality in $(4)$, [1, Theorem 3.8] gives the sharp estimate

\begin{equation*}
|\gamma_{2}|\le \frac{1}{2}(1+2e^{-2})=0.635...
\end{equation*}
\bigskip

For $n\ge 3$, the problem seems much harder, and no significant upper bounds for $|\gamma_{n}|$ when $f\in S$ appear to be known. 
\bigskip

Denote by  $S^{*}$ the subclass of $S$ of  starlike functions, so that  $f\in S^{*}$ if, and only if, for $z\in D$

\begin{equation*}
Re  \ {\dfrac{zf'(z)}{f(z)}>0.}
\end{equation*}
 \bigskip
 
\noindent Thus we can write $zf'(z)=f(z)h(z)$, where $h\in P$, the class of functions satisfying $Re \ h(z)>0$ for $z\in D$. Simple differentiation in $(1)$ again and noting that the coefficients $c_{n}$ of the Taylor series of $h$ about $z=0$ satisfy $|c_{n}|\le 2$ for $n\ge 1$, shows that $|\gamma _{n}|\le\dfrac{1}{n}$ holds for $f\in S^{*}$ and $n\ge 2$.

\bigskip
Suppose now that $f$ is analytic in $D$, then $f$ is close-to-convex if, and only if, for $z\in D$, there exists $g\in S^{*}$ such that 

\begin{equation}
Re  \ {\dfrac{zf'(z)}{g(z)}>0}.
\end{equation}
\bigskip

\noindent We denote the class of close to convex functions by $K$ and note the well-known inclusion relationship $S^{*}\subset K\subset S$.

\bigskip

That the inequality $|\gamma_{n}|\le \dfrac{1}{n}$ for $n\ge 2$ extends to the class $K$ was claimed in a paper of El Hosh [3]. However Girela [4] pointed out an error in the proof and showed that for $f\in K$, this inequality is false for $n\ge2$. In the same paper it was shown that $|\gamma_{n}|\le \dfrac{3}{2n}$ holds  for $n\ge1$ whenever $f$ belongs to  the set of the extreme points of the closed convex hull of the class K, which  implies that $|\gamma_{3}|\le \dfrac{1}{2}$ in this case.  As was pointed out above, this bound false for the entire class $K$.  It is the purpose of this paper to establish the sharp bound $|\gamma_{3}|\le \dfrac{7}{12}$ for the class $K$ when the coefficient $b_{2}$ in the Taylor expansion for $g(z)$ is real.

\bigskip

We first note that from $(4)$ it is an immediate consequence of the Fekete-Szeg\"o inequality for $f\in K$ [5] that the following sharp inequality holds for $f\in K$

\begin{equation*}
|\gamma _{2}|=\dfrac{1}{2}|a_{3}-\dfrac{1}{2}a_{2}^{2}|\le \frac{11}{18}=0.6111..
\end{equation*}
 
\bigskip

We now turn our attention to the case $n=3$ for the class $K$.

\bigskip

\noindent It follows from $(6)$ that we can write $zf'(z)=g(z)h(z)$, where $Re \ h(z)>0$ for $z\in D$ and, since $g\in S^*$, $zg'(z)=g(z)p(z)$, where $Re \ p(z)>0$ for $z\in D$.

\bigskip

\noindent Now write

\begin{equation}
h(z)=1+\sum _{n=1}^{\infty}c_{n}z^{n}
\end{equation} 
\begin{equation}
p(z)=1+\sum _{n=1}^{\infty}p_{n}z^{n}
\end{equation} 
\begin{equation}
g(z)=z+\sum_{n=2}^{\infty}b_{n}z^{n}.
\end{equation}

\bigskip

We shall need the following result [6], which has been used widely.

\bigskip

\noindent {\bf Lemma}

\bigskip

Let $h,p\in P$ and be given by $(7)$ and $(8)$ respectively, then for some complex valued $x$ with $|x|\le 1$ and some complex valued $t$ with $|t|\le 1$

\begin{equation*}
\begin{split}
2c_{2}&=c_{1}^{2}+x(4-c_{1}^{2})\\
4c_{3}&=c_{1}^{3}+2(4-c_{1}^{2})c_{1}x-c_{1}(4-c_{1}^{2})x^{2}+2(4-c_{1}^{2})(1-|x|^{2})t.\\
\end{split}
\end{equation*}

Similarly for some complex valued $y$ with $|y|\le 1$ and some complex valued $s$ with $|s|\le 1$

\begin{equation*}
\begin{split}
2p_{2}&=p_{1}^{2}+y(4-p_{1}^{2})\\
4p_{3}&=p_{1}^{3}+2(4-p_{1}^{2})p_{1}y-p_{1}(4-p_{1}^{2})y^{2}+2(4-p_{1}^{2})(1-|y|^{2})s.\\
\end{split}
\end{equation*}

We prove the following:

\bigskip

\noindent {\bf Theorem}

\bigskip

Let $f\in K$, then

$$
|\gamma _{1}|\le 1, \qquad |\gamma _{2}|\le \frac{11}{18}.
$$

Also when  $f\in K$ and  $b_{2}$ is real,

$$
\qquad |\gamma _{3}|\le \frac{7}{12}.
$$

\bigskip
The inequalities are sharp.

\bigskip
\begin{proof}

\bigskip

As noted above, the first two inequalities are proved. Thus it remains to prove the third.

\bigskip

From $(5)$ we need to find an upper bound for
\begin{equation}
|\gamma _{3}|=\frac{1}{2}|a_{4}-a_{2}a_{3}+\frac{1}{3}a_{2}^{3}|.
\end{equation}

First note that equating coefficients we have
\begin{gather*}
2a_{2}=c_{1}+p_{1}\\
3a_{3}=c_{2}+c_{1}p_{1}+\frac{p_{1}^{2}+p_{2}}{2}\\
4a_{4}=c_{3}+c_{2}p_{1}+\frac{c_{1}(p_{1}^{2}+p_{2})}{2}+\frac{p_{1}^{3}}{6}+\frac{p_{1}p_{2}}{2}+\frac{p_{3}}{3}.\\
\end{gather*}

Substituting into $(10)$ gives

\begin{equation}
\begin{split}
|a_{4}-a_{2}a_{3}+\frac{1}{3}a_{2}^{3}|&=\big|\frac{c_{3}}{4}+\frac{c_{2}p_{1}}{12}+\frac{c_{1}p_{2}}{24}+\frac{p_{3}}{12}+\frac{p_{1}p_{2}}{24}\\
&-\frac{c_{1}c_{2}}{6}-\frac{c_{1}^{2}p_{1}}{24}+\frac{c_{1}^{3}}{24}\big|.
\end{split}
\end{equation}
\bigskip

We now use the Lemma to eliminate $c_{2}$, $c_{3}$, $p_{2}$ and $p_{3}$ from $(11)$ and obtain

\begin{equation}
\begin{split}
|a_{4}-a_{2}a_{3}+\frac{1}{3}a_{2}^{3}|&=\big|\frac{c_{1}^3}{48}+\frac{c_{1}xX}{24}-\frac{c_{1}x^{2}X}{16}+\frac{XZ}{8}+\frac{p_{1}xX}{24}+\frac{c_{1}p_{1}^{2}}{48}\\
&+\frac{c_{1}yY}{48}+\frac{p_{1}yY}{16}-\frac{p_{1}y^{2}Y}{48}+\frac{p_{1}^{3}}{24}+\frac{YV}{24}\big|\\
\end{split}
\end{equation}

\noindent where for simplicity, we have set $X=4-c_{1}^{2}$, $Y=4-p_{1}^{2}$, $Z=(1-|x|^{2})s$ and $V=(1-|y|^{2})t.$

\bigskip

Without loss in generality we may write $c_{1}=c$, with $0\le c \le 2$. Also since we are assuming $b_{2}=p_{1}$ to be real, we can write $p_{1}=q$, with $0\le |q| \le 2$.  Writing $|q|=p$, it then follows using the triangle inequality in $(12)$ together with $|s|\le 1$ and $|t|\le 1$, that

\begin{equation}
\begin{split}
|a_{4}-a_{2}a_{3}+\frac{1}{3}a_{2}^{3}|&\le \frac{c^{3}}{48}+\frac{c|x|X}{24}+\frac{c|x|^{2}X}{16}+\frac{XZ}{8}+\frac{p|x|X}{24}+\frac{cp^{2}}{48}\\
&+\frac{c|y|Y}{48}+\frac{p|y|Y}{16}+\frac{p|y|^{2}Y}{48}+\frac{p^{3}}{24}+\frac{YV}{24}\\
&=F(c,p,|x|,|y|)\  \text{say.}
\end{split}
\end{equation}

\noindent where now $X=4-c^{2}$, $Y=4-p^{2}$, $Z=1-|x|^{2}$ and $V=1-|y|^{2}.$

\bigskip

Thus we need to find the maximum of $F(c,p,|x|,|y|)$ over the hyper-rectangle $R=[0,2]\times [0,2]\times[0,1]\times[0,1]$.
\bigskip

\noindent From $(13)$ substituting for $X$, $Y$, $Z$ and $V$ gives

\begin{equation}
\begin{split}
F(c,p,|x|,|y|)&=\frac{1}{48}c^{3}+\frac{1}{48}cp^{2}+\frac{1}{24}p^{3}+\frac{1}{24}c|x|(4-c^2)+\frac{1}{16}c|x|^2(4-c^2)\\
&+\frac{1}{8}(4-c^2)(1-|x|^2)+\frac{1}{24}p|x|(4-c^2)+\frac{1}{48}c|y|(4-p^2)\\
&+\frac{1}{16}p|y|(4-p^2)+\frac{1}{48}p|y|^2(4-p^2)+\frac{1}{24}(4-p^2)(1-|y|^2).
\end{split}
\end{equation}

We first assume that $F(c,p,|x|,|y|)$ has a maximum value at an interior point $(c_{0},p_{0},|x_{0}|,|y_{0}|)$ of $R$. Then since
\bigskip
\begin{equation*}
\frac{\partial F}{\partial |x|}=\frac{1}{24}c(4-c^2)+\frac{1}{4}c|x|(4-c^2)-\frac{1}{4}|x|(4-c^2)+\frac{1}{2}p(4-c^2)=0
\end{equation*}
at such a point, it follows that $c_{0}=2$, which is a contradiction. Hence any maximum points must be on the boundary of $R.$

\bigskip

\bigskip

Thus we need to find the maximum value of $F(c,p,|x|,|y|)$ on each of the $32$ edges and $24$ faces ($8$ of co-dimension $1$ and $16$ of co-dimension $2$) of $R$. Finding these maximum values involves a great many tedious exercises in elementary calculus and for the sake of brevity, we omit many of the simple ones . The process does however identify the maximum value of $7/6$ needed in the Theorem and shows that the maximum value on all edges and faces is less than or equal to $7/6$. 

\bigskip

Finding the maximum values of $F(c,p,|x|,|y|)$ on each of the $32$ edges involves trivial exercises, and shows that $F(c,p,|x|,|y|)\le 7/6$ on all of these edges. On the $16$ faces of co-dimension $2$, similar simple exercises in elementary calculus again shows that $F(c,p,|x|,|y|)\le 7/6$ on each of these faces. We thus consider the $8$ faces of co-dimension $1$  as follows.

\bigskip

On the face $c=0$, suppose that $|x|\le 1$ in $(14)$, which gives a resulting expression
\bigskip

\begin{equation*}
\begin{split}
G_{1}(0,p,|y|)&=\frac{1}{24}p^3+\frac{1}{2}+\dfrac{1}{6}p+\dfrac{1}{16}p|y|(4-p^2)\\
&+\dfrac{1}{48}p|y|^2(4-p^2)+\dfrac{1}{24}(4-p^2)(1-|y|^2).\\
\end{split}
\end{equation*}

\bigskip

\noindent Differentiating $G_{1}(0,p,|y|)$ with respect to $|y|$ shows that any maximum must occurs on the boundary of $[0,2]\times [0,1]$ and since  the largest value at the end points is $7/6$, $F(c,p,|x|,|y|)$ has maximum  $7/6$ on the face $c=0$.

\bigskip

On the face $c=2$, suppose again that $|x|\le 1$ in $(14)$, to obtain the expression
\bigskip

\begin{equation*}
\begin{split}
G_{2}(2,p,|y|)=&\frac{1}{6}+\dfrac{1}{24}p^2+\dfrac{1}{24}p^3+\dfrac{1}{24}|y|(4-p^2)\\
&+\dfrac{1}{16}p|y|(4-p^2)+\dfrac{1}{48}p|y|^2(4-p^2)+\dfrac{1}{24}(4-p^2)(1-|y|^2),
\end{split}
\end{equation*}
\noindent and following the same procedure gives a maximum of $0.696$ on $[0,2]\times [0,1]$.

\bigskip

On the face $p=0$, suppose that $|x|\le 1$ and $|y|\le 1$ in $(14)$, to obtain the expression

\begin{equation*}
\begin{split}
G_{3}(c,0,|y|)=\frac{1}{48}c^3+\dfrac{5}{48}c(4-c^2)+\dfrac{1}{8}(4-c^2)+ \dfrac{1}{12}c+\dfrac{1}{6},
\end{split}
\end{equation*}

\noindent which has maximum value $23/24$ on $[0,2]\times[0,1]$.

\bigskip

On the face $p=2$, $(14)$ becomes

\begin{equation*}
\begin{split}
G_{4}(c,2,|x|)&=\frac{1}{3}+\dfrac{c}{12}+\dfrac{1}{48}c^3+\dfrac{1}{12}(4-c^2)|x|+\dfrac{1}{24}c(4-c^2)|x|\\
&+\dfrac{1}{16}c(4-c^2)+\dfrac{1}{8}(4-c^2)(1-|x|^2).
\end{split}
\end{equation*}

\bigskip

\noindent Differentiating $G_{4}(c,2,|x|)$ with respect to $|x| $ and considering the end points gives a maximum value $1.005$ on $[0,2]\times[0,1]$.

\bigskip

On the face $|x|=0$, suppose that $|y|\le 1$ in $(14)$, to obtain

\begin{equation*}
\begin{split}
G_{5}(c,p,0)&=\dfrac{1}{48}c^3+\dfrac{1}{8}(4-c^2)+\dfrac{1}{48}cp^2+\dfrac{1}{24}p^3+\dfrac{1}{24}(4-p^2)\\
&+\dfrac{1}{12}p(4-p^2)+\dfrac{1}{48}c(4-p^2),
\end{split}
\end{equation*}

\bigskip

\noindent and it is now an easy exercise to show that $G_{5}(c,p,0)$ has a maximum value of $0.9531$ when $p=4/3$ and $c=2-2\sqrt6/3$ on $[0,2]\times[0,2]$.

\bigskip

On the face $|x|=1$, suppose that $|y|\le 1$ in $(14)$, to obtain

\bigskip

\begin{equation*}
\begin{split}
G_{6}(c,p,1)&= \dfrac{1}{48}c^3+\dfrac{1}{48}cp^2+\dfrac{1}{24}p^3+\dfrac{5}{48}c(4-c^2)+\dfrac{1}{24}p(4-c^2)\\
&+\dfrac{1}{24}(4-p^2)+\dfrac{1}{12}p(4-p^2)+\dfrac{1}{48}cp(4-p^2)\\
&=\dfrac{1}{6} + \dfrac{1}{2}c - \dfrac{1}{12}c^3 + \dfrac{1}{2}p - \dfrac{1}{24}c^2 p - \dfrac{1}{24}p^2 - \dfrac{1}{24}p^3\\
&\le \dfrac{1}{6} + \dfrac{1}{2}c - \dfrac{1}{12}c^3 + \dfrac{1}{2}p - \dfrac{1}{24}p^3.
\end{split}
\end{equation*}

\noindent It is now a simple exercise to show that this expression has maximum value $5/6$ on $[0,2]\times[0,2]$

\bigskip

On the face $|y|=0$, $(14)$ becomes

\begin{equation*}
\begin{split}
G_{7}(c,p,|x|)&=\dfrac{1}{48}c^3+\dfrac{1}{48}cp^2+\dfrac{1}{24}p^3+\dfrac{1}{24}(4-p^2)+
\\ &+\dfrac{1}{24}c|x|(4-c^2)+\dfrac{1}{24}p|x|(4-c^2)+\dfrac{1}{16}c|x|^2(4-c^2)\\
& +\dfrac{1}{8}(4-c^2)(1-|x|^2).\\
\end{split}
\end{equation*}

\bigskip

\noindent Differentiating $G_{7}(c,p,|x|)$ with respect to $|x|$ shows as before, that there are no maximum points in the interior of $[0,2]\times[0,2]\times [0,1]$, and so we need only find the maximum values of $G_{7}(c,p,|x|)$ on the boundary of $[0,2]\times[0,2]\times[0,1]$. In the interests of brevity, we omit the simple analysis which gives maximum value of $1.005$ when $p=2$ and $|x|=1$

\bigskip

We finally note that on the face $|y|=1$

\bigskip

\begin{equation*}
\begin{split}
G_{8}(c,p,|x|)&=\frac{1}{48}c^{3}+\frac{1}{48}cp^{2}+\frac{1}{24}p^{3}+\frac{1}{24}c|x|(4-c^2)+\frac{1}{16}c|x|^2(4-c^2)\\
&+\frac{1}{8}(4-c^2)(1-|x|^2)+\frac{1}{24}p|x|(4-c^2)+\frac{1}{48}c(4-p^2)\\
&+\frac{1}{12}p(4-p^2).
\end{split}
\end{equation*}

\bigskip

\noindent As before, differentiating $G_{8}(c,p,|x|)$ with respect to $|x|$ shows that there are no maximum points in the interior of $[0,2]\times[0,2]\times [0,1]$, and so we need only find the maximum values of $G_{8}(c,p,|x|)$ on the boundary of $[0,2]\times[0,2]\times[0,1]$. Again in the interests of brevity, we omit the simple analysis which gives a maximum value of $1.052$, again less than $7/6$.
\bigskip

Thus we have shown that in all cases, the maximum value of $(14)$ is at most $7/6$, which completes the proof of the Theorem.

\bigskip

We finally note that equality in the inequality in $|\gamma _{3}|\le 7/12$ is attained when $c_{1}=0$ and $c_{2}=c_{3}=p_{1}=p_{2}=p_{3}=2.$

\bigskip

\end{proof}

\bigskip

\noindent {\bf Remark 1}

\bigskip

The condition that $b_{2}$ is real in the inequality for $|\gamma _{3}|$ arises in order to maximise (12).  We conjecture that this condition can be removed and $|\gamma _{3}|\le 7/12$ for $f\in K$.
\bigskip

\noindent {\bf Remark 2}

\bigskip

The correct growth rate for $\gamma_{n}$ appears to be unknown for close-to-convex functions and in this direction the best know estimate to date appears to be that of Ye [7], who showed that $|\gamma_{n}|\le \dfrac{A\log n}{n}$, where $A$ is an absolute constant.

\end{document}